\documentclass[11pt]{amsart}
\usepackage{amssymb}
\usepackage{amsmath}
\usepackage[active]{srcltx}
\usepackage{t1enc}
\usepackage[latin2]{inputenc}
\usepackage{verbatim}
\usepackage{amsmath,amsfonts,amssymb,amsthm}
\usepackage[mathcal]{eucal}
\usepackage{enumerate}
\usepackage[centertags]{amsmath}
\usepackage{graphics}
\newtheorem{theorem}{Theorem}

\newtheorem{remark}{Remark}

\begin{document}
\author{George Tephnadze and Giorgi Tutberidze}
\title[Nörlund logaritmic means ]{A note on The maximal operators of the Nörlund logaritmic means of Vilenkin-Fourier series}
\address{George Tephnadze, The University of Georgia, School of Science and Technology, 77a Merab Kostava St, Tbilisi, 0128,  Georgia.}
\email{g.tephnadze@ug.edu.ge}
\address{G.Tutberidze, The University of Georgia, IV, 77a Merab Kostava St, Tbilisi 0128, Georgia and Department of Engineering Sciences and Mathematics, Lule\aa\ University of Technology, SE-971 87 Lule\aa, Sweden.}
\email{giorgi.tutberidze1991@gmail.com}
\thanks{The research was supported by Shota Rustaveli National Science
	Foundation grant YS-18-043.}
\date{}
\maketitle

\begin{abstract}
The main aim of this paper is to investigate $\left(H_{p},L_{p}\right)$-
type inequalities for the the maximal operators of Nörlund logaritmic means, for $0<p<1.$
\end{abstract}



\section{INTRODUCTION}

It is well-known that (see e.g. \cite{AVD}, \cite{gol} and \cite{sws}) Vilenkin systems do not form bases in the Lebesgue space $L_1\left(G_m\right).$ Moreover, there exists a function in the Hardy space $H_{1},$ such that the partial sums of $f$ are not
bounded in $L_{1}$-norm.

In \cite{tep7} (see also \cite{tepthesis}) it was proved that the following is true:

\textbf{Theorem T1.} Let $0<p<1.$ Then the maximal operator
\begin{equation*}
	\overset{\sim}{S}_{p}^{*}f :=\sup_{n\in \mathbb{N}}
	\frac{\left| S_n f \right|}{\left(n+1\right)^{1/p-1}}
\end{equation*}
is bounded from the Hardy space $H_{p}\left(G_m\right)$ to the space $L_{p}\left(G_m\right).$ Here $S_n$ denotes the $n $-th partial sum with respect to the Vilenkin system. Moreover, it was proved that the rate of the factor $(n+1)^{1/p-1}$ is in a sense sharp.

In the case $p=1$ it was proved that the maximal operator $\widetilde{S}^{\ast }$ defined by 
\begin{equation*}
	\widetilde{S}^{\ast}:=\sup_{n\in \mathbb{N}}\frac{\left\vert
		S_{n}\right\vert }{\log \left( n+1\right)}
\end{equation*}
is bounded from the Hardy space $H_{1}\left(G_m\right)$ to the space $L_{1}\left(G_m\right).$ Moreover, the rate of the factor $\log(n+1)$ is in a sense sharp. Similar problems for N\"orlund logarithmic means in the case when $p=1$ was considered in \cite{PTW3}.

M\'oricz and Siddiqi \cite{Mor} investigate the approximation properties of
some special N\"orlund means of Walsh-Fourier series of $L_{p}\left(G_m\right)$ functions in $L_p$-norm. Fridli, Manchanda and Siddiqi \cite{FMS} improved and extended the results of M\'oricz and Siddiqi \cite{Mor} to Martingale Hardy spaces. However, the
case when $\left\{q_k=1/k:k\in\mathbb{N}_+\right\}$ was excluded, since
the methods are not applicable to N\"orlund logarithmic means. In \cite{Ga2}
Gát and Goginava proved some convergence and divergence properties of
Walsh-Fourier series of the N\"orlund logarithmic means of functions in the
Lebesgue space $L_{1}\left(G_m\right).$ In particular, they proved that there exists an
function in the space $L_1\left(G_m\right),$ such that 
\begin{equation*}
\sup_{n\in \mathbb{N}}\left\Vert L_{n}f\right\Vert _{1}=\infty.
\end{equation*}

In \cite{bg1} (see also \cite{PTW3,tep10}) it was proved that there exists a martingale $f\in H_{p}\left(G_m\right), \quad (0<p<1)$ such that 
\begin{equation*}
\sup_{n\in \mathbb{N}}\left\| L_{n}f\right\|_p=\infty.
\end{equation*}

Analogical problems for N\"orlund means with respect to Walsh, Kaczmarz and unbounded Vilenkin systems were considered in \cite{GNCz}, \cite{bt1, bt2}, \cite{NT1, NT2, NT3}, \cite{PTW1, PTW2}, \cite{tep8, tep9,tepthesis}, \cite{tut1}.

In this paper we discuss boundedness of weighted maximal
operators from the Hardy space $H_p\left(G_m\right)$ to the Lebesgue space $L_p\left(G_m\right)$ for $0<p<1.$

\section{Definitions and Notation}

Let $\mathbb{N}_{+}$ denote the set of the positive integers, $\mathbb{N}:=\mathbb{N}_{+}\cup \{0\}.$

Let $m:=(m_{0},m_{1},...)$ denote a sequence of the positive integers not
less than $2.$

Denote by
\begin{equation*}
Z_{m_{k}}:=\{0,1,...m_{k}-1\}
\end{equation*}
the additive group of integers modulo $m_{k}.$

Define the group $G_{m}$ as the complete direct product of the group $%
Z_{m_{j}}$ with the product of the discrete topologies of $Z_{m_{j}}$ $^{,}$%
s.

The direct product $\mu $ of the measures
\begin{equation*}
\mu _{k}\left( \{j\}\right) :=1/m_{k}\text{ \qquad }(j\in Z_{m_{k}})
\end{equation*}
is the Haar measure on $G_{m_{\text{ }}}$with $\mu \left( G_{m}\right) =1.$

If $\sup\limits_{n\in \mathbb{N}}m_{n}<\infty $, then we call $G_{m}$ a bounded Vilenkin
group. If the generating sequence $m$ is not bounded then $G_{m}$ is said to
be an unbounded one. \textbf{In this paper we discuss bounded
Vilenkin groups only.}

The elements of $G_{m}$ are represented by sequences
\begin{equation*}
x:=(x_0,x_1,...x_j,...)\qquad \left(x_k\in Z_{m_{k}}\right).
\end{equation*}

It is easy to give a base for the neighborhood of $G_{m}$
\begin{equation*}
I_{0}\left( x\right) :=G_{m},
\end{equation*}
\begin{equation*}
I_{n}(x):=\{y\in G_{m}\mid y_{0}=x_{0},...y_{n-1}=x_{n-1}\}\text{ }(x\in
G_{m},\text{ }n\in \mathbb{N})
\end{equation*}
Denote $I_n:=I_n\left(0\right),$ for $n\in\mathbb{N}$ and $\overset{-}{I_n}:=G_{m}$ $\backslash $ $I_n$ .

If we define the so-called generalized number system based on $m$ in the
following way :
\begin{equation*}
M_{0}:=1,\text{ \qquad } M_{k+1}:=m_{k}M_{k} \qquad (k\in \mathbb{N})
\end{equation*}
then every $n\in \mathbb{N}$ can be uniquely expressed as $n=\overset{\infty }{
\underset{k=0}{\sum }}n_{j}M_{j}$ where $n_{j}\in Z_{m_{j}}$ $(j\in \mathbb{N})$ and
only a finite number of $n_{j}`$s differ from zero. Let $\left| n\right|
:=\max\{j\in \mathbb{N}; \quad n_j\neq 0\}.$

The norm (or quasi-norm) of the space $L_{p}(G_{m})$ is defined by 
\begin{equation*}
\left\| f\right\| _{p}:=\left( \int_{G_{m}}\left| f\right| ^{p}d\mu\right) ^{1/p}\qquad \left( 0<p<\infty \right) .
\end{equation*}
The space $weak-L_{p}\left( G_{m}\right) $ consists of all measurable
functions $f$ for which
\begin{equation*}
\left\| f\right\| _{weak-L_{p}(G_{m})}:=\underset{\lambda >0}{\sup }\lambda
^{p}\mu \left( x:\left| f\left( x\right) \right| >\lambda \right) <+\infty.
\end{equation*}

Next, we introduce on $G_{m}$ an ortonormal system which is called the
Vilenkin system.

At first define the complex valued function $r_{k}\left( x\right)
:G_{m}\rightarrow C,$ the generalized Rademacher functions as
\begin{equation*}
r_{k}\left(x\right):=\exp\left( 2\pi ix_{k}/m_{k}\right) \qquad
\left( i^{2}=-1,\qquad x\in G_{m},\qquad k\in \mathbb{N}\right).
\end{equation*}

Now define the Vilenkin system $\psi :=(\psi _{n}:n\in \mathbb{N})$ on $G_{m}$ as:
\begin{equation*}
\psi_{n}:=\overset{\infty }{\underset{k=0}{\Pi }}r_{k}^{n_{k}}, \qquad \left( n\in \mathbb{N}\right).
\end{equation*}

Specifically, we call this system the Walsh-Paley one if m=2.

The Vilenkin system is ortonormal and complete in $L_{2}\left( G_{m}\right)
\,$\cite{AVD,Vi}.

Now we introduce analogues of the usual definitions in Fourier-analysis.

If $f\in L_{1}\left( G_{m}\right) $ we can establish the the Fourier
coefficients, the partial sums of the Fourier series, the Dirichlet kernels
with respect to the Vilenkin sistem $\psi $ in the usual manner:
\begin{eqnarray*}
\widehat{f}&:&=\int_{G_{m}}f \overline{\psi }_{k}d\mu,\qquad \left(k\in \mathbb{N}\right), \\
S_{n}f&:&
=\sum_{k=0}^{n-1}\widehat{f}\left(k\right)\psi_{k}, \qquad \left(n\in \mathbb{N}_{+}, \quad S_{0}f:=0\right), \\
D_{n}&:&=\sum_{k=0}^{n-1}\psi _{k},\qquad
\left(n\in \mathbb{N}_{+}\right).
\end{eqnarray*}

Recall that (for details see e.g. \cite{AVD})
\begin{equation} \label{3}
D_{M_{n}}\left( x\right) =\left\{
\begin{array}{l}
M_{n}\quad x\in I_{n} \\
0 \qquad x\notin I_{n}
\end{array}
\right. 
\end{equation}

The $\sigma$-algebra generated by the intervals 
$\left\{ I_n\left(x\right) :x\in G_{m}\right\} $ will be denoted by $\digamma _{n}$ $\left(n\in \mathbb{N}\right).$ Denote by $f=\left( f_n:n\in \mathbb{N}\right) $ a
martingale with respect to $\digamma _{n}$ $\left( n\in \mathbb{N}\right).$ (for details see e.g. \cite{We1, We3}). The maximal function of a martingale $f$ is defend by 
\begin{equation*}
f^{*}=\sup_{n\in \mathbb{N}}\left| f_n\right| .
\end{equation*}

In the case when $f\in L_{1},$ the maximal function is also be given by
\begin{equation*}
f^{*}\left( x\right)=\sup_{n\in \mathbb{N}}\frac{1}{\left| I_{n}\left( x\right)\right| }\left| \int_{I_{n}\left( x\right) }f\left( u\right) \mu \left(u\right) \right|.
\end{equation*}

For $0<p<\infty $ the Hardy martingale spaces $H_{p}$ $\left( G_{m}\right) $
consist of all martingales for which
\begin{equation*}
\left\| f \right\|_{H_p}:=\left\| f^{*}\right\|_{p}<\infty.
\end{equation*}

If $f\in L_{1},$ then it is easy to show that the sequence 
$\left(S_{M_n}f :n\in \mathbb{N}\right) $ is a martingale. If $f=\left(f_n:n\in \mathbb{N}\right) $ is martingale then the Vilenkin-Fourier coefficients must be defined in a slightly different manner: 
\begin{equation*}
\widehat{f}\left(i\right):=\lim_{k\rightarrow \infty} \int_{G_{m}}f_k\overline{\psi}_id\mu.
\end{equation*}

The Vilenkin-Fourier coefficients of $f\in L_{1}\left( G_{m}\right) $ are
the same as those of the martingale $\left( S_{M_{n}}f :n\in
\mathbb{N} \right) $ obtained from $f$ .

Let \{$q_{k}: k>0$\} be a sequence of nonnegative numbers. The $n$-th 
Nörlund means for the Fourier series of $f$ is defined by
\begin{equation*}
\frac{1}{Q_n}\overset{n}{\underset{k=1}{\sum }}q_{n-k}S_k f
\qquad \text{where} \qquad
Q_{n}:=\overset{n}{\underset{k=1}{\sum }}q_{k}.
\end{equation*}

If $q_{k}=1/k,$ then we get Nörlund logarithmic means
\begin{equation*}
L_{n}f:=\frac{1}{l_{n}}\overset{n-1}{
\underset{k=0}{\sum }}\frac{S_{k}f}{n-k} \qquad \text{where} \qquad l_{n}=\overset{n-1}{\underset{k=0}{\sum }}\frac{1}{n-k}=\overset{n}{%
\underset{j=1}{\sum }}\frac{1}{j}.
\end{equation*}

A bounded measurable function $a$ is $p$-atom, if there exist a dyadic
interval $I,$ such that 
\begin{equation*}
\int_{I}ad\mu=0, \ \  
\left\|a\right\| _{\infty }\leq \mu \left(
I\right) ^{-1/p}, \ \  
\text{supp}\left(a\right) \subset I.
\end{equation*}

\section{Formulation of Main Results}

\begin{theorem}
a) Let $0<p<1.$ Then the maximal operator
\begin{equation*}
\overset{\sim}{L}_{p}^{*}f :=\sup_{n\in \mathbb{N}}
\frac{\left| L_n f \right|}{\left(n+1\right)^{1/p-1}}
\end{equation*}
is bounded from the Hardy space $H_{p}\left( G_{m}\right) $ to the space $%
L_{p}\left( G_{m}\right).$

b) Let $0<p<1$ and $\varphi :\mathbb{N}_{+}\rightarrow [1, \infty)$ be a
non-decreasing function satisfying the condition
\begin{equation} \label{6}
\overline{\lim_{n\rightarrow\infty}}\frac{n^{1/p-1}}{\log n\varphi\left( n\right)}=+\infty.  
\end{equation}
Then there exists a martingale $f\in H_p\left(G_m\right),$ such that the maximal operator
\begin{equation*}
\text{ }\sup_{n\in \mathbb{N}}\frac{\left| L_{n}f 
\right| }{\varphi \left( n+1\right) }
\end{equation*}
is not bounded from the Hardy space $H_{p}\left( G_{m}\right) $ to the space
$L_{p}\left( G_{m}\right).$
\end{theorem}

\section{Proof of the Theorem}

\begin{proof} Since 
\begin{equation*}
\frac{\left|L_{n} f \right|}{\left(n+1\right)^{1/p-1}}\leq \frac{1}{\left(n+1\right)^{1/p-1}}\underset{1\leq k\leq n}{\sup }\left|S_k f  \right|\leq \underset{1\leq k\leq n}{\sup }\frac{\left|S_k f \right|}{\left(k+1\right)^{1/p-1}}\leq \underset{n\in \mathbb{N}}{\sup }\frac{\left|S_n f  \right|}{\left(n+1\right)^{1/p-1}}
\end{equation*}
if we use Theorem T1 we obtain that 
\begin{equation*}
\underset{n\in \mathbb{N}}{\sup }\frac{\left|L_n f \right|}{\left(n+1\right)^{1/p-1}}
\leq \underset{n\in \mathbb{N}}{\sup }\frac{\left|S_n f  \right|}{\left(n+1\right)^{1/p-1}}
\end{equation*}
and 
\begin{equation*}
\left\Vert \underset{n\in \mathbb{N}}{\sup }\frac{\left|L_{n} f \right|}{\left(n+1\right)^{1/p-1}}\right\Vert_{p}\leq \left\Vert \underset{n\in \mathbb{N}}{\sup }\frac{\left|S_{n} f \right|}{\left(n+1\right)^{1/p-1}}\right\Vert _{p}\leq c_{p}\left\Vert f\right\Vert_{H_p}.
\end{equation*}
Now, prove part b) of the Theorem. Let
\begin{equation*}
f_{n_k}=D_{M_{2n_k+1}}-D_{M_{2n_k}}.
\end{equation*}

It is evident
\begin{equation}
\widehat{f}_{n_{k}}\left( i\right) =\left\{
	\begin{array}{l}
		\text{ }1,\qquad \text{ if }\qquad i=M_{2n_{k}},...,M_{2n_{k}+1}-1, \\
		\text{ }0,\qquad \text{otherwise.}
	\end{array}
	\right.  \label{14a}
\end{equation}
Then we can write
\begin{equation} \label{14}
S_{i}f_{n_{k}} =\left\{
\begin{array}{l}
D_{i}-D_{M_{2n_{k}}}, \quad \text{ if } \quad
i=M_{2n_{k}}+1,...,M_{2n_{k}+1}-1, \\
f_{n_{k}},\qquad \qquad \quad \text{ if } \quad i\geq M_{2n_{k}+1}, \\
0,\qquad \qquad \qquad \text{ otherwise.}
\end{array}
\right.  
\end{equation}

From (\ref{3}) we get
\begin{eqnarray} \label{14c}
\left\| f_{n_{k}}\right\| _{H_{p}}&=&\left\|
\sup\limits_{n\in \mathbb{N}}S_{M_{n}}f_{n_{k}}
\right\|_{p}  =\left\| D_{M_{2n_{k}+1}}-D_{M_{_{2n_k}}} \right\|_p \\
&\leq& \left\| D_{M_{2n_{k}+1}} \right\|_p+ \left\|D_{M_{_{2n_{k}}}}\right\| _p\leq cM_{_{2n_k}}^{1-1/p}<c<\infty.  \notag
\end{eqnarray}

Let $0<p<1$ and $\left\{ \lambda_k:k\in \mathbb{N}_{+}\right\}$ be an increasing sequence of the positive integers such that
\begin{equation*}
\lim_{k\rightarrow\infty}\frac{\lambda_{k}^{1/p-1}}{\varphi \left(
\lambda _{k}\right)}=\infty.
\end{equation*}

Let$\ \left\{ n_{k}:k\in \mathbb{N}_{+}\right\} \subset \left\{ \lambda_k:k\in
\mathbb{N}_{+}\right\} $ such that
\begin{eqnarray*}
\underset{k\rightarrow\infty}{\lim}\frac{\left(
M_{_{2n_{k}}}+2\right)^{1/p-1}}{\log{\left( M_{2n_{k}}+2\right)}\varphi\left(
M_{2n_k+2}\right)}\geq c\lim_{k\rightarrow \infty }\frac{\lambda _{k}^{1/p-1}}{\varphi \left(\lambda _{k}\right) }=\infty.
\end{eqnarray*}

According to (\ref{14}) we can conclude that
\begin{eqnarray*}
&&\left| \frac{L_{M_{2n_{k}}+2}f_{n_{k}}}{\varphi\left(
M_{2n_k+2}\right)} \right|=\frac{\left| D_{M_{_{2n_{k}}}+1} -D_{M_{_{2n_{k}}}}\right|}{l_{M_{2n_{k}}+1}\varphi\left(M_{2n_k+1}\right)}\\
&=&\frac{\left|\psi _{M_{_{2n_{k}}}} \right|}{
l_{M_{2n_{k}}+2}\varphi\left(
M_{2n_k+1}\right)}=\frac{1}{l_{M_{2n_{k}}+1}\varphi\left(
M_{2n_k+2}\right)}.
\end{eqnarray*}

Hence,
\begin{equation} \label{17}
\mu \left\{ x\in G_{m}:\left| L_{M_{2n_{k}}+2}f_{n_{k}}\right|\geq \frac{1}{l_{M_{2n_{k}}+2}\varphi\left(
	M_{2n_k+2}\right)}\right\}  
=\mu \left( G_{m}\right) =1.
\end{equation}

By combining (\ref{14c}) and (\ref{17}) we get that

\begin{eqnarray*}
&&\frac{\frac{1}{l_{M_{2n_{k}}+2}\varphi\left(M_{2n_k+2}\right)} \left( \mu \left\{ x\in G_{m}:\quad\left| L_{M_{2n_{k}}+2} f_{n_{k}}  \right| 
\geq\frac{1}{l_{M_{2n_{k}}+2}\varphi\left(
M_{2n_k+2}\right)}\right\}\right)^{1/p}}{\left\| f_{n_{k}}\right\| _p} \\
&\geq &\frac{M_{_{2n_{k}}}^{1/p-1}}{l_{M_{2n_{k}}+2}\varphi\left(
M_{2n_k+2}\right)} \geq \frac{c\left(M_{_{2n_{k}}}+2\right)^{1/p-1}} {\log{\left(M_{2n_{k}}+2\right)}\varphi\left(M_{2n_k+2}\right)}\rightarrow \infty, \text{\quad as \quad }k\rightarrow \infty.
\end{eqnarray*}

Theorem is proved.

\end{proof}

\textbf{Open Problem.} For any  $0<p<1$ let find non-decreasing function $\Theta :\mathbb{N}_{+}\rightarrow [1, \infty)$ suth that the following maximal operator

\begin{equation*}
\overset{\sim}{L}_{p}^{*}f :=\sup_{n\in \mathbb{N}}
\frac{\left| L_n f \right|}{\Theta\left(n+1\right)}
\end{equation*}
is bounded from the Hardy space $H_{p}\left( G_{m}\right) $ to the Lebesgue space $L_p\left(G_m\right)$ and the rate of $\Theta :\mathbb{N}_{+}\rightarrow [1,\infty)$ is sharp, that is, for any non-decreasing function
$\varphi :\mathbb{N}_{+}\rightarrow [1, \infty)$ satisfying the condition
\begin{equation*} 
\overline{\lim_{n\rightarrow \infty }}\frac{\Theta\left( n\right)}{ \varphi\left( n\right)}=+\infty,  
\end{equation*}
then there exists a martingale $f\in H_p\left( G_m\right),$ such that the maximal operator
\begin{equation*}
\sup_{n\in\mathbb{N}}\frac{\left|L_nf\right|}{\varphi\left(n+1\right)}
\end{equation*}
is not bounded from the Hardy space $H_p\left(G_m\right)$ to the space
$L_p\left(G_m\right).$

\begin{remark}
According to Theorem 1 we can conclude that there exist absolute constants $C_1$ and $C_2$ such that
\begin{equation*}
\frac{C_1 n^{1/p-1}}{\log(n+1)}\leq\Theta\left(n\right) \leq C_2 n^{1/p-1}.
\end{equation*}
\end{remark}

\end{document}